# Normal Extensions of an Singular Differential Operator for First Order

## by

## Elgiz BAIRAMOV[1)], Rukiye ÖZTÜRK MERT[2)] & Zameddin ISMAILOV[2)]

[1)] Ankara University, Faculty of Sciences, Department of Mathematics, 06100, Ankara, TURKEY

[2)] Karadeniz Technical University, Faculty of Sciences, Department of Mathematics, 61080, Trabzon, TURKEY

e-mail: bairamov@science.ankara.edu.tr; rukiye-55@hotmail.com; zameddin@yahoo.com

**Abstract**

In this work, in the Hilbert space of vector-functions $L^2(H, (-\infty, a) \cup (b, +\infty)), a < b$ all normal extensions of the minimal operator generated by linear singular formally normal differential expression $l(\cdot) = (\frac{d}{dt} + A_1, \frac{d}{dt} + A_2)$ with a selfadjoint operator coefficients $A_1$ and $A_2$ in any Hilbert space $H$, are described in terms of boundary values. Later structure of the spectrum of these extensions is investigated.

## 1. İntroduction

Many problems arising in the modelling of processes of multi-particle quantum mechanics, quantum field theory, in the physics of rigid bodies, in the theory of multipoint boundary value problems for differential equations and ets. support to study normal extension of formally normal differential operators in direct sum of Hilbert spaces([1-3]).

The general theory of normal extensions of formally normal operators in any Hilbert space and their spectral theory have been investigated by many mathematicians (see, for example [4-9] ). Applications of this theory to two point differential operators in Hilbert space of functions are continied today even($[10-17]$).

İt is clear that, the minimal operators $L_{10}$ and $L_{20}$ generated by differential expressions in forms

$$\frac{d}{dt} + A_k, k = 1,2$$

in the Hilbert spaces of vector-functions $L_2(H, (-\infty, a)), L_2(H, (b, \infty))$, respectively, where $A_1 = A_1^* \leq 0, A_2 = A_2^* \geq 0, -\infty < a < b < \infty$ ,are maximal formally normal, consequently have not normal extensions.

Note that the table is changed in the multipoint case in the next sense. Direct sum $L_{10} \oplus L_{20}$ of operators $L_{10}$ and $L_{20}$ is formally normal ,is not maximal, moreover it have normal extensions in the direct sum $L_2(H,(-\infty,a)) \oplus L_2(H,(b,\infty))$.

In this paper will be considered a linear multipoint differential expression

$$l = (l_1, l_2), l_k = \frac{d}{dt} + A_k, k = 1,2$$

in the direct sum of Hilbert spaces of vector-functions

$$L_2 = L_2(H,(-\infty,a)) \oplus L_2(b,\infty))$$

where, $A_1 = A_1^* \leq 0, A_2 = A_2^* \geq 0, -\infty < a < b < \infty$.

The main purpose of this work is to describe in terms of boundary values all normal extensions of the minimal operator $L_0 = L_{10} \oplus L_{20}$ and investigate structure of the spectrum of such extensions. Note that in the regular case, i.e. in the finite interval description and spectral analysis of normal extensions have been established in [10].

## 2. Description of Normal Extensions.

Let $H$ be a separable Hilbert space and $a, b \epsilon R, a < b$. In the Hilbert space $L_2 = L_2(H,(-\infty,a)) \oplus L_2(H,b,\infty))$ of vector-functions consider the following linear multipoint differential expression in form

$$l(u) = (l_1(u_1), l(u_2)) = (du_1/dt + A_1 u_1, du_2/dt + A_2 u_2), \quad u = (u_1, u_2) \quad (2.1)$$

where $A_k: D(A_k): H \to H, k = 1,2$ are linear selfadjoint operators and $A_1 \leq 0, A_2 \geq 0$.

In this case the formally adjoint differntial expression be in form

$$l^+(v) = (l_1^+(v_1), l_2^+(v_2)) = (-dv_1/dt + A_1 v_1, -dv_2/dt + A_2 v_2), \quad v = (v_1, v_2) \quad (2.2)$$

Define by $L_{k0}(L_{k0}^+), k = 1,2$ the minimal operator generated by differential expression(2.1)((2.2)) and $L_k = (L_{k0}^+)^*(L_k^+ = L_{k0}^*) k = 1,2$ the maximal operator generated by differential expression(2.1)((2.2)) in the space $L_2(H,(-\infty,a))$ ($L_2(H,b,\infty))$).In this case the operators $L_0 = L_{10} \oplus L_{20}$ and $L = L_1 \oplus L_2$ in the space $L_2$ are called a minimal and maximal operators generated by (2.1) and (2.2), respectively. It is clear that $L_{k0} \subset L_k$, $L_{k0}^+ \subset L_k^+, k = 1,2$ and $L_0 \subset L$.

**Definition 2.1.** A densely defined closed operator $N$ in any Hilbert space is called formally normal if $D(N) \subset D(N^*)$ and $\|Nf\| = \|N^*f\|$ for all $f \epsilon D(N)$. If a formally normal operator has no formally normal extension, then it is called maximal formally normal operator. If a formally normal operator $N$ satisfies the condition $D(N) = D(N^*)$, then it is called normal operator ([4]).

**Definition 2.2.** Let $S$ be a linear closed symmetric operator in the Hilbert space $H$ with finite or infinite deficiency indices. Then triplet $(\mathcal{H}, Y_1, Y_2)$, where $\mathcal{H}$ is a Hilbert space and $Y_1$ and $Y_2$ are linear mapping from $D(S^*)$ to $\mathcal{H}$ is called a space of boundary values of the operator $S$:

(1) for any $u, v \in D(S^*)$,

$$(S^*u, v)_H - (u, S^*v)_H = (Y_1(u), Y_2(v))_{\mathcal{H}} - (Y_2(u), Y_1(v))_{\mathcal{H}};$$

(2) for any $g, h \epsilon \mathcal{H}$ there exists element $f \epsilon D(S^*)$ such that

$$Y_2(f) = g, \ Y_2(f) = h.$$

Note that every closed symmetric operator in any Hilbert space with equal deficiency indices has at least one a space of boundary values (see[18]).

Now construct a space of boundary values for the minimal operator $M_0$ generated by linear singular differential expression of first order in form

$$m(u) = (m(u_1), m(u_2)) = \left(-i\frac{du_1}{dt}, -i\frac{du_2}{dt}\right)$$

in the direct sum $L_2(H, (-\infty, a)) \oplus L_2(H, (b, \infty))$. Note that the minimal operator $M_0$ is closed symmetric operator in $L_2(H, (-\infty, a)) \oplus L_2(H, (b, \infty))$ with deficiency indices $(dimH, dimH)$ which from it will be used in last investigation.

**Theorem. 2.3.** The triplet $(H, Y_1, Y_2)$, where

$$Y_1: D(M_0^*) \to H, Y_1(u) = \frac{1}{i\sqrt{2}}(u_2(b) + u_1(a)),$$

$$Y_2: D(M_0^*) \to H, Y_2(u) = \frac{1}{\sqrt{2}}(u_2(b) - u_1(a)), \quad u = (u_1, u_2) \in D(M_0^*)$$

İs a space of boundary values of the minimal operator $M_0$ in $L_2$.

**Proof.** For arbitrary $u = (u_1, u_2)$ and $v = (v_1, v_2)$ from $D(M_0^*)$ validity the following equality

$$(M_0^* u, v)_{L_2} - (u, M_0^* v)_{L_2} = (Y_1(u), Y_2(v))_H - (Y_2(u), Y_1(v))_H$$

can be easily verified. Now give any elements $f, g \in H$. Find the function $u = (u_1, u_2) \in D(M_0^*)$ such that

$$Y_1(u) = \frac{1}{i\sqrt{2}}(u_2(b) + u_1(a)) = f \quad \text{and} \quad Y_2(u) = \frac{1}{\sqrt{2}}(u_2(b) - u_1(a)) = g,$$

that is,

$$u_1(a) = (if - g)/\sqrt{2} \quad \text{and} \quad u_2(b) = (if + g)/\sqrt{2}.$$

If choose these functions $u_1(t), u_2(t)$ in following form

$$u_1(t) = \int_{-\infty}^{t} e^{s-a}\, ds (if - g)/\sqrt{2}, \quad t < a, \quad u_2(t) = \int_{t}^{\infty} e^{b-t}\, ds\, (if + g)/\sqrt{2}, \quad t > b,$$

then it is clear that $(u_1, u_2) \in D(M_0^*)$ and $Y_1(u) = f$, $Y_2(u) = g$. ∎

Firstly we prove the following.

**Theorem 2.4.** If the minimal operators $L_{10}, L_{20}$ are formally normal, then are true

$$D(L_{10}) \subset W_2^1(H,(-\infty,a)), A_1 D(L_{10}) \subset L_2(H,(-\infty,a))$$

$$D(L_{20}) \subset W_2^1(H,(b,\infty)), A_2 D(L_{20}) \subset L_2(H,(b,\infty))$$

**Proof.** Indeed, in this case for each $u_1 \in D(L_{10}) \subset D(L_{10}^*)$ is true

$u_1' + A_1 u_1 \in L_2(H,(-\infty,a))$ and $u_1' - A_1 u_1 \in L_2(H,(-\infty,a))$, hence

$u_1' \in L_2(H,(-\infty,a))$ and $A_1 u_1 \in L_2(H,(-\infty,a))$, i.e.,

$D(L_{10}) \subset W_2^1(H,(-\infty,a))$ and $A_1 D(L_{10}) \subset L_2(H,(-\infty,a))$. ∎

In similar way it is proved the second part of theorem.

**Theorem 2.5.** Let $(-A_1)^{1/2} W_2^1(H,(-\infty;a)) \subset W_2^1(H,(-\infty,a))$ and $A_2^{1/2} W_2^1(H,(b,\infty))$
$\subset W_2^1(H,(b,\infty))$. Each normal extention $\tilde{L}$ of the minimal operator $L_0$ in the Hilbert space $L_2$ is generated by differential expression (2.1) and boundary contition

$$u_2(b) = W u_1(a), u_1(a) \in \ker(-A_1)^{1/2}, u_2(b) \in \ker A_2^{1/2}, \quad (2.3)$$

where $W: H \to H$ is a unitary operator. Moreover, the unitary operator $W$ in $H$ is determined by the extension $\tilde{L}$, i.e. $\tilde{L} = L_W$ and vice versa.

**Proof.** Let $\tilde{L}$ be a normai extension of $L_0$. Then the operator $Im\tilde{L}$ is selfadjoint extension of the minimal operator $Im L_0$ in $L_2$. Consequently, operator $Im\tilde{L}$ is generated by the differential expression $-id/dt$ and the boundary condition

$$(W - E)Y_1(u) + i(W + E)Y_2(u) = 0, u = (u_1,u_2) \in D(L),$$

where $W: H \to H$ is a unitary operator and is unique $(see[18])$. From this later simple calculations it is obtained that

$$u_2(b) = W u_1(a) \text{ for all } u = (u_1, u_2) \in D(\tilde{L}), \text{ i.e. } \tilde{L} = L_W.$$

On the other hand since extension $\tilde{L}$ is normal, then from the equality $\|\tilde{L} u\| = \|\tilde{L}^* u\|$ for every $u = (u_1, u_2) \in D(\tilde{L})$ we have

$$(u', A_1 u)_{L_2(H,(-\infty,a))} + (A_1 u, u')_{L_2(H,(-\infty,a))} + (u', A_2 u)_{L_2(H,(b,\infty))} + (A_2 u, u')_{L_2(H,(b,\infty))} = 0$$

and from this

$$\|(-A_1)^{1/2} u(a)\|^2 + \| A_2^{1/2} u(b) \|^2 = 0$$

Hence it is obtained that $u(a) \in ker(-A_1)^{1/2}, u(b) \in kerA_2^{1/2}$.

On the contrary, let $L_W$ a extension of the minimal operator $L_0$ which is generates by linear expression $l(.)$ and boundary condition (2.3).Then it is easily to consider that a adjoint operator $L_W^*$ which is generated by differential expression $l^+(.)$ and boundary condition

$$v_1(a) = W^*v_2(b), v = (v_1, v_2) \in D(L^+)$$

is a normal operator in $L_2$. In this case can be easily verify that an operator $L_W$ is extension of the $L_0$ and it is normal. ∎

**Corollary 2.6.** If at least one of the operators $A_1$ and $A_2$ is one-to-one mapping in $H$, then minimal operator $L_0$ is maximally formal normal in $L_2$.

**Corollary 2.7.** If there exists at least one normal extension of the minimal operator $L_0$, then

$$dimker(-A_1)^{1/2} = dimkerA_2^{1/2} > 0.$$

### 3. THE SPECTRUM of the NORMAL EXTENSİONS

In this section structure of the spectrum of the normal extension $L_W$ in $L_2$ will be investigated. Later will be assumed that $A_1 = A_1^* \leq 0, A_2 = A_2^* \geq 0$ and $0 \in \sigma_p((-A_1)^{1/2}) \cap \sigma_p(A_2^{1/2})$.

Firstly prove the following proposition.

**THEOREM 3.1:** The point spectrum of the any normal extension $L_W$ of the minimal operator $L_0$ in $L_2$ is empty, i.e.

$$\sigma_p(L_W) = \emptyset.$$

**Proof.** Consider the following problem for the normal extension $L_W$

$$L_W u = \lambda u, \lambda = \lambda_r + i\lambda_i \in \mathbb{C}, u \in D(L_W).$$

From this we have

$$u_1' + A_1 u_1 = \lambda u_1, u_1 \in L_2(H, (-\infty, a)),$$

$$u_2' + A_2 u_2 = \lambda u_2, u_2 \in L_2(H, (b, \infty)), \lambda \in \mathbb{C}.$$

$$u_2(b) = Wu_1(a), u_1(a) \in ker(-A_1)^{1/2}, u_2(b) \in kerA_2^{1/2}.$$

The general solutions of last equations in form

$$u_1(\lambda; t) = e^{-(A_1 - \lambda)(t-a)} f_\lambda, t < a, f_\lambda \in H_{-1/2}((-A_1))$$

$$u_2(\lambda; t) = e^{-(A_2 - \lambda)(t-b)} g_\lambda, t > b, g_\lambda \in H_{-1/2}(A_2),$$

$$g_\lambda = Wf_\lambda, f_\lambda = u_1(\lambda; a), g_\lambda = u_2(\lambda; b).$$

Since $0 \in \sigma_p((-A_1)^{1/2}) \cap \sigma_p(A_2^{1/2})$ and $(-A_1)^{1/2} f_\lambda = 0$, $A_2^{1/2} g_\lambda = 0$, then we have

$$u_1(\lambda; t) = e^{\lambda(t-a)} f_\lambda, t < a, f_\lambda \in H_{-1/2}((-A_1))$$

$$u_2(\lambda; t) = e^{\lambda(t-b)} g_\lambda, t > b, g_\lambda \in H_{-1/2}(A_2),$$

$$g_\lambda = W f_\lambda, \ f_\lambda = u_1(\lambda; a), g_\lambda = u_2(\lambda; b).$$

In order to $u_1(\lambda; .) \in L_2(H, (-\infty, a))$ and $u_2(\lambda; .) \in L_2(H, (b, \infty))$ necessary and sufficient condition is $\lambda_r \geq 0$ and $\lambda_r \leq 0$ respectively. Hence $\lambda_r = 0$.

Consequently,

$$u_1(\lambda; t) = e^{i\lambda_i(t-a)} f_\lambda, t < a,$$

$$u_2(\lambda; t) = e^{i\lambda_i(t-b)} g_\lambda, t > b, g_\lambda = W f_\lambda.$$

In this case it is clear that for the $u_1(\lambda; .) \in L_2(H, (-\infty, a))$ and $u_2(\lambda; .) \in L_2(H, (b, \infty))$ necessary and sufficient conditions are $f_\lambda = 0, g_\lambda = 0$. From this implies that $u_1 = 0$ and $u_2 = 0$ in $L^2$. Therefore $\sigma_p(L_W) = \emptyset$. ∎

Since resudial spectrum of any normal operator in any Hilbert space is empty, then furthermore the continuous spectrum of normal extensions $L_W$ of the minimal operator $L_0$ is investigated.

**Theorem 3.2.** For the continuous spectrum of the normal extension $L_W$ is true
$$\sigma_c(L_W) = iR.$$

**Proof.** Assume that $\lambda \in \sigma_c(L_W)$. Then on the important theorem for the spectrum of normal operators [19], that is,

$$\sigma(L_W) \subset \sigma(ReL_W) + i\sigma(ImL_W),$$

it is obtained that

$$\lambda_r \in \sigma(ReL_W), \lambda_i \in \sigma(ImL_W).$$

From this imply that $\lambda_r \in \sigma(A_1), \lambda_r \in \sigma(A_2)$, hence on the conditions to the operators $A_1$ and $A_2$ we have $\lambda_r = 0$. On the other hand from the proof of previes theorem it is to see that $\ker(L_W - \lambda) = \{0\}$ for any $\lambda \in \mathbb{C}$. Consequently, $\sigma_c(L_W) \subset iR$.

Furthermore, it is clear that for the $\lambda = i\lambda_i \in \mathbb{C}$ the general solution of the boundary value problem

$$u_1' + A_1 u_1 = i\lambda_i u_1 + f_1, \quad u_1, f_1 \in L_2(H, (-\infty, a)),$$

$$u_2' + A_2 u_2 = i\lambda_i u_2 + f_2, \quad u_2, f_2 \in L_2(H, (b, \infty)), \lambda_i \in R,$$

$$u_2(b) = W u_1(a), u_1(a) \in \ker(-A_1)^{1/2}, u_2(b) \in \ker A_2^{1/2}$$

Is in form

$$u_1(i\lambda_i; t) = e^{-(A_1-i\lambda_i)(t-a)} f_{i\lambda_i} - \int_t^a e^{-(A_1-i\lambda_i)(t-s)} f_1(s)ds, \quad t < a,$$

$$u_2(i\lambda_i; t) = e^{-(A_2-i\lambda_i)(t-b)} g_{i\lambda_i} + \int_b^t e^{-(A_2-i\lambda_i)(t-s)} f_2(s)ds, \quad t > b,$$

$$g_{i\lambda_i} = W f_{i\lambda_i}.$$

In this case it is true

$e^{-(A_1-i\lambda_i)(t-a)} f_{i\lambda_i} \in L_2(H, (-\infty, a))$, $e^{-(A_2-i\lambda_i)(t-b)} g_{i\lambda_i} \in L_2(H, (b, \infty))$ for any $g_{i\lambda_i}, f_{i\lambda_i} \in H$.

Here if choose $f_1(t) = e^{i\lambda_i t} e^{-(t-a)} f^*, f^* \in \ker(-A_1)^{1/2}, t < a$, then

$$\int_t^a e^{-(A_1-i\lambda_i)(t-s)} f_1(s)ds = e^{-i\lambda_i t} \int_t^a e^{-(s-a)} f^* ds = e^{-i\lambda_i t}(e^{-(t-a)} - 1)f^*, t < a.$$

Therefore

$$\int_{-\infty}^a \|e^{-i\lambda_i t}(e^{-(t-a)} - 1)f^*\|^2 dt = \int_{-\infty}^a \|e^{-i\lambda_i t}(e^{-(t-a)} - 1)f^*\|^2 dt$$

$$= \int_{-\infty}^a (e^{-2(t-a)} - 2e^{-(t-a)} + 1)dt \|f^*\|^2 = \infty$$

Consequently, for the such $f_1(t) \in L_2(H, (-\infty, a))$, $u_1(i\lambda_i; t) \notin L_2(H, (-\infty, a))$. This means that for any $\lambda \in \mathbb{C}$, $a$ $operator$ $L - \lambda$ is one-to-one in $L_2$, but it is not onto transformation. On the other hand, since resudial spectrum $\sigma_r(L_W)$ is empty, then it is implies that $\sigma(L_W) = \sigma_c(L_W) = iR$. ∎

**Example 3.3.** By the last theorem the spectrum of following boundary value problem

$$\partial u(t,x)/\partial t - \text{sgnt}\, (\partial^2 u(t,x))/(\partial x^2) = f(t,x), |t|>1, x \in [0,1],$$

$$u(1,x) = e^{i\varphi} u(-1,x), \varphi \in [0, 2\pi),$$

$$\partial u(t,0)/\partial t = \partial u(t,1)/\partial t = 0, |t|>1,$$

in the space $L_2((-\infty,-1) \times (0,1)) \oplus L_2((1,\infty) \times (0,1))$ is continuous and coincides with $iR$.